\documentclass{nyj}
\usepackage{amsmath,amsfonts,epsfig}
\usepackage{amsthm}
\usepackage{amssymb}
\newtheorem{prop}{Proposition}
\newtheorem{cor}[prop]{Corollary}

\newtheorem{lem}[prop]{Lemma}

\newtheorem{theo}[prop]{Theorem}

\theoremstyle{definition}
\newtheorem{defn}[prop]{Definition}

\newcommand{\longversion}[1]{}
\newcommand{\eg}{{\underline EG}}
\newcommand{\gcw}{{$G$-$CW$}}
\newcommand{\NN}{\mathbb{N}}

\newcommand{\ra}{\rightarrow}

\title[$\eg$ for Hyperbolic Groups]{A Model for the Universal Space
  for Proper Actions of a Hyperbolic 
  Group} 
\author{David Meintrup}
\address{Universit{\"a}t der Bundeswehr, M{\"u}nchen, Germany}
\email{david.meintrup@unibw-muenchen.de}

\author{Thomas Schick} 
\address{Universit{\"a}t G{\"o}ttingen, Germany}
\email{schick@uni-math.gwdg.de, www.uni-math.gwdg.de/schick/}

\date{2001, December 19}
\thanks{The first author has been
  supported by the DAAD.}

\keywords{universal space for proper actions, Rips complex, word
  hyperbolic group, Gromow hyperbolic group, fixed point set,
  conjugacy classes of finite subgroups, finiteness properties for
  universal spaces, classifying space for proper actions, G-CW-complex}

\subjclass{20F67, 55R35, 57M07}

\begin{document}

\begin{abstract}{
Let $G$ be a word hyperbolic group in the sense of Gromov and $P$ its
associated Rips complex. We prove that the fixed point set $P^H$ is
contractible for every finite subgroups $H$ of $G$. This is the main ingredient
for proving that $P$ is a finite model for the universal space $\eg$ of proper
actions. As a corollary we get that a hyperbolic group has only finitely many 
conjugacy classes of finite subgroups.}

\end{abstract}

\maketitle

\section{Introduction}

The aim of this note is to prove the following theorem.

\begin{theo} \label{Rips is eg-underbar}
Let $G$ be a $\delta$-hyperbolic group relative to a set of generators~$S$, 
$P_d(G,S)$ its Rips complex, with $d \geq 32\delta + 20$,
and $P$ its second barycentric subdivision. Then
$P$ is a finite \gcw-model for the universal space $\eg$ of proper actions, 
i.e.~$P$ consists of
finitely many $G$-equivariant cells, $P$ has only finite
stabilizers and $P^H$ is contractible for any finite subgroup $H
\subset G$.
\end{theo}

This has been claimed in \cite[section 2]{Baum-Connes-Higson(1994)},
but there does not seem to be a proof in the literature. We will give
a complete and
detailed proof. We are grateful for fruitful discussions with
Martin Bridson. We thank the referee for some useful suggestions about
the presentation.

We start with some basic definitions about universal spaces for proper
actions. For a more detailed introduction see
\cite{Lueck-Meintrup(2000)}. Let $G$ be a discrete group. A
\gcw-complex $X$ is a $CW$-complex with a $G$-action that is cellular and
whenever a cell is mapped on itself by some $g \in G$, the restriction of the
action of $g$ to 
this cell is the identity map. We call a \gcw-complex $X$ \emph{proper} if
the stabilizers $G_x:=\{g\in G\mid gx=x\}$ are finite for all $x \in X$. A \gcw-complex is called 
\emph{finite} if $X/G$ is compact. \longversion{(I.e. $X$ consists of
  finitely many equivariant cells.)}

\begin{defn}
A \emph{\gcw-model for the universal space $\eg$ for proper actions of $G$}
is a proper \gcw-complex such that the fixed point sets $\eg^H$ are
contractible for every finite subgroup $H$ of $G$.
\end{defn}

For every group $G$ there is such a model $\eg$ and any two such models are
$G$-homotopy equivalent. The space $\eg$ has the following universal
property: For any proper \gcw-complex $X$, an up to $G$-homotopy
unique $G$-map $X \ra \eg$ exists. 
\longversion{In other words, $\eg$ is the final object in the 
$G$-homotopy category of proper \gcw-complexes.} We mention that the space
$\eg$ plays an important role in the formulation of the Baum-Connes
Conjecture \cite{Baum-Connes-Higson(1994)} and for the generalization of the 
Atiyah-Segal completion theorem \cite{Lueck-Oliver(1998)}. In particular,
  finiteness conditions for $\eg$ have been studied for discrete groups, e.g. in
  \cite{Kropholler-Mislin(1998)}, \cite{Lueck(2000)}, 
  and for locally compact groups in \cite{Lueck-Meintrup(2000)}.

\section{Hyperbolicity and the Rips Theorem}

Let $G$ be a finitely generated group and $S$ a finite
set of generators. We
will always assume that the identity element is not contained in $S$,
 $e \notin S$  and that $S$ is symmetric, that is $S = S^{-1}$.
Given $x,y\in G$, write $x^{-1}y=s_1\cdots s_d$ with a minimal number $d$ 
of $s_i\in S$. Then $d_S(x,y):=d$ gibes the left-invariant word metric 
on $G$. We set $[x,y]:=\{x, xs_1,xs_1s_2,\ldots,xs_1\cdots
s_d=y\}$. This is a ``geodesic'' joining $x$ and $y$. Note that
$[x,y]$ is not unique.
See \cite[1.2]{Ghys-delaHarpe(1990)} for more details.

\begin{defn}
Let $d \in \NN$. The {\it Rips complex } $P_d(G, S)$ is the simplicial
complex whose $k$-simplices are given by $(k+1)$-tuples
 $(g_0,\ldots, g_k)$ of pairwise distinct elements of $G$
with $\max\limits_{0 \leq i,j \leq k}d_S(g_i,g_j) \leq d$. 
\longversion{We equipe $P_d(G, S)$
with the weak topology.} Observe that the $0$-skeleton of $P_d(G,S)$
coincides with $G$.
\end{defn}

If no confusion is possible, we well omit the notation of the set of
generators, and simply write $d(.,.)$ instead of $d_S(.,.)$ for the word
metric and $P_d(G)$ instead of $P_d(G,S)$ for the Rips complex.
For finite subsets $K,L$ of vertices we will use the notation
$$
d(K,L) := \max_{k \in K,l \in L}d(k,l)
$$
for the maximal distance between $K$ and $L$. The {\it diameter} of $K$ is
given by 
$$
d(K) := d(K,K) = \max_{k,k' \in K} d(k,k').
$$
Typically, we will look at diameters of orbits $Hx$ of finite subgroups
$H$. Notice that the $H$-invariance of the word metric implies the identity
$d(Hx) = d(\{x\},Hx)$.
Since the word metric is left invariant, we have a simplicial action
of $G $ on $P_d(G, S)$ given by
$g\cdot (g_0,\ldots,g_k)
= (g g_o,\ldots,g g_k).$

\begin{defn}
Let $X$ be a metric space and $\delta \geq 0$. Then $X$ is
 {\it $\delta$-hyperbolic} if for any four points  $x, y, z, t \in X$ the following
 inequality holds:
\[
d(x,y) + d(z,t) \leq \max\{d(x,z) + d(y,t), d(x,t)+d(y,z)\} + 2 \delta.
\]
A group $G$ is $\delta$-hyperbolic if it is $\delta$-hyperbolic as metric
 space equipped with the word metric.
\end{defn}

We quote the following theorem from
\cite[ 12.- Th{\'e}or{\`e}me, p. 73]{Ghys-delaHarpe(1990)}:

\begin{theo}[Rips-Theorem] \label{Rips}
Let $G$ be a $\delta$-hyperbolic group for the set of generators $S$.
Let $d \geq 4\delta+2$ and $P$ be the second barycentric subdivision of
$P_d(G, S)$. Then $P$ is a contractible, locally finite simplicial
complex of finite dimension and with a simplicial action of $G$ which is
faithful and properly discontinuous. Moreover:
\begin{itemize}
\item[(i)] The stabilizer of each simplex is finite.
\item[(ii)] If $p$ is a vertex of $P$ and $g \in G$
such that $p\not= g p$, then the stars of $p$ and $g p $
 are disjoint.
\item[(iii)] The orbit space $G \backslash P$ is a finite
simplicial complex
and the projection $\pi: P \rightarrow G \backslash P $ is simplicial.
\item[(iv)] If $G$ is torsionfree then the action is free.
\end{itemize}
\end{theo}

The Rips-Theorem already proves some parts of Theorem \ref{Rips
is eg-underbar}. Also, we can conclude from parts (i) and (iv), that $P$ is
a finite 
\gcw-model for $\eg= EG$ in the torsion-free case. For the general case,    
the next lemma is an important tool, because it provides the existence
of universally small orbits. In the following, for a real number $r$ the
notation $[r]$ will always mean the integral part of $r$.  

\begin{lem} \label{small orbits}
Let $G$ be a $\delta$-hyperbolic group and $y_0$  a vertex of the Rips
complex $P_d(G)$. Let $H$ be a finite subgroup and $d(Hy_0) = R$. 
\begin{itemize}
\item[(a)] Then there is a vertex $x$ of $P_d(G)$ such that
\begin{eqnarray*}
d(Hx,Hy_0) & \leq & \left \lfloor\frac{R}{2}\right \rfloor + 2 \delta + 1 \\
\mbox{and } d(Hx) & \leq & 8 \delta + 4.
\end{eqnarray*}
\item[(b)] If, in addition $R \geq 8\delta + 2 $ and $d(y_0,x_0) =
d(Hy_0, x_0)$ for some vertex $x_0$ of $P_d(G)$, then
$$
d(Hx,x_0) \leq d(x_0, y_0).
$$
\end{itemize}
\end{lem}
\begin{proof}
Let $y' \in Hy_0$ be a vertex on the orbit such that
$d(y',y_0) = R$ and let $x$ be a vertex on a geodesic $[y',y_0]$ such that
$d(x,y_0) = \left [\frac{R}{2}\right ]$.  Then, using the definition of hyperbolicity for the
points $x,hy_0,y_0,y'$ we get:
\begin{eqnarray*}
d(x,hy_0) & \leq & \max\{d(y',hy_0) + d(x,y_0), d(y_0,hy_0) +
d(x,y')\} - d(y_0,y') + 2 \delta \\
& \leq & R + \left \lfloor\frac{R}{2}\right \rfloor + 1 - R + 2 \delta \\
& \leq & \left \lfloor\frac{R}{2}\right \rfloor + 2\delta + 1. 
\end{eqnarray*}
Now the left invariance of the metric gives the first statement of
$(a)$. For the second part let $r:= d(x,hx) $ for an $h \in H$ and 
$v$ be a vertex on a geodesic $[x,hx]$ such that $d(v,x) =
\left [\frac{r}{2}\right ]$. Since $d(Hy_0) = R$, there is a vertex 
$z \in Hy_0$ such that
$d(v,z) \geq \left [\frac{R}{2}\right ]$. Applying the hyperbolicity to the
points $x,hx,v,z$ we get:
\begin{eqnarray*}
r  = d(x,hx) & \leq & \max\{d(x,z)+ d(v,hx), d(hx,z) + d(v,x)\} - d(v,z) +
2\delta \\
& \leq & \left \lfloor \frac{R}{2}\right \rfloor + 2 \delta + 1 + \frac{r}{2} + 1 - \left \lfloor\frac{R}{2}\right \rfloor + 2 \delta \\
& = & \frac{r}{2} + 4 \delta + 2. 
\end{eqnarray*}
This implies $r \leq 8 \delta + 4$ what we wanted to show. 

For part $(b)$ of the lemma we apply hyperbolicity to $hx,x_0,y_0,y'$ and
get
\begin{eqnarray*}
d(hx,x_0) & \leq & \max\{d(hx,y_0) + d(x_0,y'), d(hx,y') + d(x_0,y_0) \}
-d(y_0,y') + 2 \delta\\
& \leq & \left \lfloor\frac{R}{2}\right\rfloor + 2\delta + 1 + d(x_0,y_0) - R + 2 \delta \\
& \leq & d(x_0,y_0)
\end{eqnarray*}
because we assumed $R \geq 8\delta + 2$. 
\end{proof}

The next proposition shows the contractibility of the fixed point sets.
It is, roughly speaking, an $H$-invariant version of the proof of the
non-equivariant
contractibility of the Rips complex, cf. \cite[4.2]{Ghys-delaHarpe(1990)}.
\begin{prop} \label{contractible fixed points}
Let $G$ be $\delta$-hyperbolic with set of generators $S$ and
$P_d(G)$ its Rips complex with $d \geq 32\delta + 20$. Let $H$ be a finite
subgroup of $G$. Then $P_d(G)^H$ is contractible. 
\end{prop}

\begin{proof}
Let $F$ be the following subcomplex of $P_d(G)$. A simplex $\sigma$ of
$P_d(G)$ belongs to $F$ if $\sigma$ contains an $H$-fixed point. This is
the case if and only if $H$ permutes the vertices of $\sigma$. Now add all
the faces of these simplices to make $F$ a subcomplex. Clearly, $F$ is an
$H$-invariant subcomplex and $F^H = P_d(G)^H$. To show the contractibility
of $F^H$ we proceed as follows: Let $K'$ be a finite subcomplex of $F$,
and set $K := HK'$.
We will show that $K$ is $H$-equivariantly contractible in $F$, 
i.e. the inclusion $K \hookrightarrow F$ is $H$-equivariantly
homotopic to a constant map. This implies that 
$\pi_i(F^H) = 0 $ for $i\geq 0$. Since $F^H$ is simplicial (after two
barycentric subdivisions), this is all we
need to prove. 

By Lemma \ref{small orbits} we can find a vertex $x_0$ with $d(Hx_0) \leq 8
 \delta + 4$. Hence $x_0 \in F$, i.e. $F$ is not empty.
 Without loss of generality we can
 assume that $x_0 \in K$. Let $K_0$ denote the finite set
 of vertices of $K$.  We distinguish two cases:
\begin{itemize}
\item[(i)] $\max_{y \in K_0} d(x_0,y) \leq \frac{d}{2}$. Then $K_0$ spans an
 $H$-invariant simplex $\sigma$ that contains a fixed point. 
Any $H$-equivariant
contraction of $\sigma$ to this fixed point will also contract $K$ in $\sigma$. 
\item[(ii)] $\max_{y \in K_0} d(x_0,y) > \frac{d}{2}$.
Let $y_0 \in K_0$ be a point furthest away from $x_0$, i.e. 
$$
d(y_0,x_0) = \max_{y \in K_0} d(x_0,y).
$$

Let $y_0'$ be the point on a geodesic $[x_0,y_0]$ with
$$
d(y_0',x_0) = d(y_0,x_0) - \left \lfloor\frac{d}{4}\right \rfloor.
$$
Next, we want to define the function
\begin{eqnarray*}
f_0: (K_0,x_0) \rightarrow (F,x_0), \quad  f(hy_0) & := & hy_0', \ h \in H, \\
  f(y) & := & y, \ y \in K_0\backslash Hy_0.
\end{eqnarray*}
Notice that $H$ acts freely on $K_0$, hence the first part of the
definition of $f$ makes sense. 
We have to verify two more things to know that this function is
well-defined. First that
$x_0 \notin Hy_0$. This follows from the fact that $d(Hx_0) \leq 8
\delta + 4 \leq \frac{d}{2}$, but $d(x_0,y_0) > \frac{d}{2}$. 
Secondly, that
$y_0' \in F$. We will show this by proving $d(Hy_0') \leq d$. Then, by
definition of $F$, we have $y_0' \in F$. Again, we have to look at two
cases:
\begin{itemize}
\item[(a)] $d(Hy_0) \leq \frac{d}{2}$: Then we have by the triangle inequality
\begin{eqnarray*}
d(hy_0',y_0') & \leq & d(hy_0',hy_0) + d(hy_0,y_0) + d(y_0,y_0') \\
 & \leq & \left \lfloor\frac{d}{4}\right \rfloor + \frac{d}{2} + \left \lfloor\frac{d}{4}\right \rfloor \leq d,
\end{eqnarray*}
hence $d(Hy_0') \leq d$. 
\item[(b)] $d(Hy_0) > \frac{d}{2}$: Since $y_0\in F$ we know that
  $d(Hy_0) \leq d$. Thus we can apply
Lemma \ref{small orbits}(a) to $Hy_0$ to obtain a vertex $x$ with orbit
$Hx$ satisfying 
\begin{equation*}
d(Hx) \leq 8 \delta +
4
\end{equation*}
and
$$
d(Hx,y_0) \leq \frac{d}{2} + 2 \delta + 1.
$$
Since we  have $d(Hy_0) >  \frac{d}{2} \geq 8 \delta + 4$ and $d(x_0,y_0) =
d(x_0, Hy_0) $ by the choice of $y_0$, Lemma \ref{small orbits}(b) also gives 
$$
d(Hx,x_0) \leq d(x_0,y_0).
$$
Hence, applying hyperbolicity to the points $hx,y_0',y_0,x_0$ we get:
\begin{eqnarray*}
d(hx,y_0') & \leq &  \max \left \{ 
\begin{array}{l} d(hx,y_0) + d(y_0',x_0)-d(y_0,x_0), \\
{d(hx,x_0)+d(y_0,y_0') - d(y_0,x_0)}
\end{array} 
\right \} + 2 \delta \\ 
& \leq & \max \{\frac{d}{2} +2 \delta + 1 - \left \lfloor\frac{d}{4}\right \rfloor, \left \lfloor\frac{d}{4}\right \rfloor \} + 2 \delta \\
& \leq & \frac{d}{2} - \left \lfloor\frac{d}{4}\right \rfloor + 4 \delta + 1 \\
& \leq & \frac{d}{2}.
\end{eqnarray*}
The last estimation holds because we assumed $d \geq 32 \delta + 20$, hence
$\left [\frac{d}{4}\right ] \geq 4 \delta + 1$. Now we can again use the triangle inequality:
$$
d(hy_0',y_0') \leq d(hy_0',hx) + d(hx,y_0') \leq \frac{d}{2} + \frac{d}{2} = d.
$$
and we get $d(Hy_0') \leq d$. 
\end{itemize}
Now that we know that $f_0$ is well defined we claim:  

{\bf Claim:} $f_0$ can be extended to a simplicial map $f:(K,x_0) \rightarrow (F,x_0)$. 

We have to show that for $x,y \in K_0$, $d(x,y) \leq d$ implies
$d(f_0(x),f_0(y)) \leq d$. 
Since we only moved the orbit of $y_0$, there is only one non-trivial case
to check,
the implication:
\begin{equation} \label{last}
d(y,hy_0) \leq d \Rightarrow d(y,hy_0') \leq d.
\end{equation}
Because of the left-invariance of the metric (replace $y$ by $h^{-1}y$) this is equivalent to:
$$
d(y,y_0) \leq d \Rightarrow d(y,y_0') \leq d.
$$
 Applying hyperbolicity to $y,y_0',y_0,x_0$ we get
\begin{eqnarray*}
d(y,y_0') & \leq  & \max \left \{ \begin{array}{l} 
 d(y,y_0) + d(x_0,y_0') -d(y_0,x_0), \\ 
d(y_0',y_0) + d(y,x_0) - d(y_0,x_0) 
\end{array}
\right \} + 2 \delta \\
& \leq & \max \{ d - \left \lfloor\frac{d}{4}\right \rfloor, \left \lfloor\frac{d}{4}\right \rfloor \} + 2 \delta \\
& = & d - \left \lfloor\frac{d}{4}\right \rfloor + 2 \delta \leq d,
\end{eqnarray*}
since $d \geq 32 \delta + 20$ implies $\left \lfloor\frac{d}{4}\right \rfloor \geq 2 \delta$.

Now $f$ is by definition an $H$-equivariant map. It remains to show that
$f$ is $H$-homotopic to the inclusion map. But this follows by
noticing that for any simplex $\sigma$ of $K$ the set 
$f(\sigma) \cup \sigma$ is contained in a simplex of $F$. This is clear
except for the case where one point is in the orbit of $y_0'$, but then it
follows from implication (\ref{last}).

Finally for each $h\in H$ we have
\begin{equation*}
d(f(hy_0),x_0)=d(f(y_0),h^{-1}x_0)\le
d(y_0',x_0)+d(x_0,h^{-1}x_0)\le d(y_0,x_0)-\left \lfloor\frac{d}{4}\right
\rfloor+8\delta+4<d,
\end{equation*}
since $d\ge 32\delta+20$, so that $\left\lfloor
  \frac{d}{4}\right\rfloor >8\delta+4$. Therefore, the whole orbit of
$y_0$ is moved closer to $x_0$.

Continuing this process on the finite complex $f(K)$ and iterating finitely
many times leads to a finite subcomplex that will satisfy
case~(i). This ends the proof. 
\end{itemize}
\end{proof} 

\begin{proof}[Proof of Theorem \ref{Rips is eg-underbar}] 
Since $P$ is a second barycentric subdivision of a simplicial complex with
a simplicial $G$-action, $P$ is by \cite[p.~117]{Bredon(1972)} a \gcw-complex. 
Theorem \ref{Rips}(i)
states that all stabilizers are finite, (iii) proves that $P$ is
finite as \gcw-complex. The contractibility of $P^H$ for finite groups
$H$ is shown in Proposition \ref{contractible fixed points}.
\end{proof}

The next corollary generalizes \cite[Prop. 4.13, p. 73]{Ghys-delaHarpe(1990)}.
An alternative proof can be found in \cite[Theorem 3.2]{Bridson-Haefliger(1999)}.
\begin{cor} \label{number_of_con}
A $\delta$-hyperbolic group $G$ has a finite number of conjugacy classes of
finite subgroups.
\end{cor}
\begin{proof} This is true for every discrete group  with a model for $\eg$ of finite
type, as is shown in \cite[Theorem 4.2]{Lueck(2000)}. 
\end{proof}

%\bibliography{lit}

\begin{thebibliography}{GdlH90}

\bibitem[BCH94]{Baum-Connes-Higson(1994)}
Paul Baum, Alain Connes, and Nigel Higson.
\newblock Classifying space for proper actions and ${K}$-theory of group
  ${C}\sp \ast$-algebras.
\newblock In {\em $C\sp \ast$-algebras: 1943--1993 (San Antonio, TX, 1993)},
  pages 240--291. Amer. Math. Soc., Providence, RI, 1994.

\bibitem[BH99]{Bridson-Haefliger(1999)}
Martin~R. Bridson and Andr{\'e} Haefliger.
\newblock {\em Metric spaces of non-positive curvature}.
\newblock Springer-Verlag, Berlin, 1999.

\bibitem[Bre72]{Bredon(1972)}
Glen~E. Bredon.
\newblock {\em Introduction to compact transformation groups}.
\newblock Academic Press, New York, 1972.
\newblock Pure and Applied Mathematics, Vol. 46.

\bibitem[GdlH90]{Ghys-delaHarpe(1990)}
Etienne Ghys and Pierre de~la Harpe, editors.
\newblock {\em Sur les groupes hyperboliques d'apr{\`e}s {M}ikhael {G}romov}.
\newblock Birkh{\"a}user Boston Inc., Boston, MA, 1990.

\bibitem[KM98]{Kropholler-Mislin(1998)}
Peter~H. Kropholler and Guido Mislin.
\newblock Groups acting on finite-dimensional spaces with finite stabilizers.
\newblock {\em Comment. Math. Helv.}, 73(1):122--136, 1998.

\bibitem[LM00]{Lueck-Meintrup(2000)}
Wolfgang L{\"u}ck and David Meintrup.
\newblock On the universal space for groups acting with compact isotropy.
\newblock In {\em Proceedings of the Conference on Geometry and Topology,
  {A}arhus, 1998}, pages 293--306. AMS Proceedings, 2000.

\bibitem[LO98]{Lueck-Oliver(1998)}
Wolfgang L{\"u}ck and Robert Oliver.
\newblock The completion theorem in ${K}$-theory for proper actions of a
  discrete group.
\newblock Preprintreihe {SFB} 478 -- {G}eometrische {S}trukturen in der
  {M}athematik, {M}{\"u}nster, Heft 1, 1998.

\bibitem[L{\"u}c00]{Lueck(2000)}
Wolfgang L{\"u}ck.
\newblock The type of the classifying space for a family of subgroups.
\newblock {\em J. Pure Appl. Algebra}, 149(2):177--203, 2000.

\end{thebibliography}

\end{document}